Bernard Ycart[*]

# 1827 : la mode de la statistique en France

origine, extension, personnages

**Résumé.** De manière largement indépendante du développement scientifique de la discipline, une mode de la statistique s'est développée en France, à partir de 1827. L'événement déclencheur a probablement été la « Carte figurative de l'instruction populaire » de Charles Dupin, avec sa fameuse ligne Saint-Malo Genève, censée séparer le nord instruit du sud ignorant. Produire, sous le nom de « statistique », des descriptions plus ou moins chiffrées sur des sujets les plus variés, est devenu un vecteur d'ambition, un moyen privilégié d'accès à la notoriété. Au-delà des traces littéraires, le phénomène peut être mesuré par sa pénétration sémantique dans la presse. Alors que l'ambition de la plupart des statisticiens amateurs restait strictement descriptive, certains ont su poser la question de la démonstration par les chiffres. C'est d'autant plus remarquable que, au sein de la science institutionnelle, les techniques de démonstration statistique, introduites par Laplace dès la fin du siècle précédent, ont longtemps été largement ignorées.

**Abstract. 1827: the trend of statistics in France; origin, extension, characters.** Independent to a great extent from the scientific development of the discipline, a trend for statistics has developed in France, from 1827 on. It was probably sparked by Charles Dupin's 'Carte figurative de l'instruction populaire', with its famous Saint-Malo Geneva line, supposed to separate the educated North from the ignorant South. It became attractive to produce, under the name 'statistics', more or less quantitative descriptions on any subject. Beyond literary records, the phenomenon can be measured by its semantic penetration in the press. Even if the ambition of most of these amateurs has remained strictly descriptive, some of them did raise the issue of proving through numbers. This is particularly remarkable, since within institutional science, the techniques of statistical proving, that had been introduced by Laplace at the end of the 18th century, have remained largely ignored for a very long time.

[*] Laboratoire Jean Kuntzmann, Université Grenoble Alpes, 51 rue des mathématiques, 38041 Grenoble cedex. E-mail: bernard.ycart@imag.fr

Passé l'« âge d'or de la statistique régionale[1] » sous Napoléon, la Restauration avait amené à la discipline, son lot de difficultés politiques et de mutations idéologiques[2]. Même pour un grand savant comme Fourier[3], il n'était plus de bon ton sous Louis XVIII, de s'adonner à la statistique. Au début du règne de Charles X, quelle ambition pouvaient avoir les héritiers de ces érudits de province qui avaient si activement contribué aux statistiques de leurs départements sous le consulat[4] ? Combien parmi eux auraient parié sur l'essor de la discipline ? Lesquels auraient pu deviner qu'ils préféreraient bientôt se déclarer « statisticien » plutôt que savant, ou bien qu'ils publieraient leurs observations archéologiques ou botaniques sous le nom de « statistiques ». Pourtant, en 1830, Balzac dans un article intitulé « De la mode en littérature » écrivait :

> Enfin, le moindre cacographe est membre d'une société savante, et ceux qui ne savent rien ou ne peuvent pas écrire comptent les fontaines de Paris, examinent les couleurs des numéros que le préfet impose aux maisons, et se prétendent occupés de statistique ; car la statistique est devenue à la mode, et c'est une position que de statistiquer[5].

L'histoire de cette mode, son origine, son développement, les traces qu'elle a laissées dans la presse et la littérature, ainsi que son rapport avec le développement de la discipline scientifique, sont l'objet de cet article.

Le phénomène est décelable dès 1827. Il a précédé, accompagné et sans doute facilité l'« enthousiasme statistique[6] » des institutions à partir de la révolution de Juillet, mais en est largement indépendant. Sur le plan épistémologique, inutile d'y chercher un nouveau paradigme : il ne s'agit le plus souvent que de rebaptiser « statistiques » des compte-rendus descriptifs, éventuellement non chiffrés. Sur le plan sociologique, le personnage du notable de province, dont l'ambition sociale et la soif de reconnaissance s'exprime au travers de la statistique, devient très vite un archétype, reconnu au point de servir de cible à la satire dans les romans ou les pièces de théâtre, dès 1827 et jusqu'à la fin du XIXe. Le déclencheur a probablement été la fameuse ligne Saint-Malo Genève séparant, selon Charles Dupin (1784–1873)[7], le nord « éclairé » du sud ignorant. Les réactions à sa fameuse « carte ombrée[8] », tout particulièrement dans les départements « noircis », ne se sont pas limitées à des protestations indignées. De nombreux notables de province ont voulu répondre à Dupin sur son propre terrain, celui de la statistique. Les premières années, le département des Bouches-du-Rhône est à la pointe du mouvement. Il fournit en 1829 la première livraison de la « Statistique morale de la France[9] », en fait une biographie départementale, conçue par Antoine Andraud (1795–1859). On y retrouve la base sociologique des « gloires locales », qui modernise la figure du philosophe, savant ou académicien provincial du siècle précédent[10].

Au delà des témoignages contemporains dans la presse, et des citations littéraires attestant du phénomène, une quantification a été tentée. Le Journal des Débats Politiques et Littéraires[11] entre 1814 et 1848, a été choisi comme corpus de référence. L'analyse textuelle fait apparaître une brusque augmentation de la proportion des numéros dans lesquels apparaît le terme « statistique », à partir de 1827 : de une à deux occurrences par mois environ jusqu'en 1826, on passe à une à deux par semaine à partir de 1828. Diverses interprétations de cette observation seront proposées. Il semble néanmoins possible de conclure que l'augmentation de fréquence du terme dans le journal des débats ne faisait qu'accompagner sa plus grande notoriété parmi le lectorat, liée à un accroissement des publications utilisant le terme dans leur titre ou leur sous-titre. Une analyse des occurrences dans la presse et dans la littérature, permettra de montrer l'extension du phénomène, ainsi que les réticences qu'il a pu susciter.

---

1. L'expression est de Perrot, J.-C., 1976, qui reste une référence indispensable sur la statistique napoléonienne, l'autre étant Bourguet, M.-N., 1988.
2. Analysées en particulier dans Brian, É., 1991, et 1994, chapitre IV.
3. Arago, F., 1854, p. 359.
4. Perrot, J.-C., 1976, p. 255.
5. La Mode, 29 mai 1830, reproduit par Therenty, M.-È., 2014. Cité par Bru B. & M.-F. & Bienaymé, O., 1997, p. 174.
6. Rosanvallon, P., 1995, p. 210.
7. Les deux références importantes sur Charles Dupin sont Christen, C. & Vatin, F., 2009, et Bradley, M. 2012. Sur la ligne Saint-Malo Genève, voir dans la première référence la contribution de Pierre Karila-Cohen, p. 129.
8. Parmi les nombreuses références sur la fameuse carte choroplète de Dupin, voir Lepetit, B., 1986, Palsky, G. 2008, de Falguerolles, A., 2010, Boria, E., 2013.
9. Andraud, A., 1829.
10. Roche, D., 1978.
11. Désigné dans tout ce qui suit par Journal des Débats, accessible sur gallica.bnf.fr.

Nombre de critiques, y compris dans la presse, posent la question de fond : que prouvent vraiment les statistiques ? La théorie des probabilités avait été développée par Pierre-Simon Laplace (1749–1827), plusieurs décennies auparavant. Dès ses premiers mémoires, Laplace avait posé les bases de la statistique inférentielle : il avait décrit le moyen de calculer des niveaux de confiance, des p-valeurs[12] ; bref, quiconque avait compris Laplace connaissait le moyen de rendre rigoureuse une affirmation basée sur des chiffres, c'est-à-dire de quantifier une incertitude ou d'évaluer une significativité. Pourtant, il aura fallu pratiquement un siècle après Laplace pour que la statistique devienne scientifique. Non pas que Laplace ait manqué de porte-paroles au sein même de la science officielle : les historiens n'ont pas manqué de souligner les rôles de Poisson, Fourier, Bienaymé, Cournot et autres[13], mais aussi leur impuissance à imposer la rigueur statistique. Comment d'ailleurs douter de l'« évidence mathématique », quand l'exemple du manque de rigueur venait du sommet, à savoir de celui qui depuis son « coup médiatique » de 1826 était l'incarnation même du statisticien que l'on n'osait contredire : Charles Dupin. Dans ces conditions, comment un statisticien amateur, loin des cercles institutionnels de la science parisienne, sans le bagage mathématique qui lui aurait permis de comprendre la difficile théorie de Laplace, aurait-il pu se sentir tenu à une méthodologie rigoureuse tant dans le recueil des données que dans leur interprétation ? Sans aller jusqu'à des développements mathématiques qui étaient hors de leur portée, certains amateurs ont pourtant laissé le témoignage de leur honnêteté intellectuelle. Les travaux de deux d'entre eux seront présentés : Jules Honoré Modeste Bigeon (1803–1831)[14], et Adolphe d'Angeville (1796–1856)[15]. Le second en particulier, développe une technique originale à base de statistiques de rangs, un siècle avant que les techniques non paramétriques deviennent d'actualité[16]. Même si, faute d'outils mathématiques, il ne va pas au bout de sa logique, ses affirmations sont bel et bien démontrables statistiquement.

En conclusion, une évaluation des étendues géographique, sociologique et historique du phénomène sera proposée.

# 1. La carte de Dupin et ses conséquences

## *La ligne Saint-Malo Genève*

L'histoire de la « Carte figurative de l'instruction populaire », que Charles Dupin présente au Conservatoire des Arts et Métiers le mercredi 29 novembre 1826, est bien connue[17] ; nous n'y revenons que pour en préciser l'écho. La carte et le discours du Conservatoire sont reproduits dans

---

12  Sur le développement par Laplace de la théorie des probabilités et son application à la statistique inférentielle, voir Stigler, S., 1986, Bru, B. 1986, Hald, A. 1998, Kuusela, V., 2012.
13  Entre autres Bru, B. & M.-F. & Bienaymé, O. 1997, Heyde, C. & Seneta, E. 1977.
14  Bigeon, J. H. M., 1829.
15  D'Angeville, A., 1837.
16  Droesbeke J.-J. & Fine, J. 1996, chapitre I.
17  Voir Palsky, G. 2008 et Boria, E., 2013 pour le point de vue des géographes, Lepetit, B., 1986 pour celui des économistes, de Falguerolles A., 2010 pour celui des statisticiens, parmi bien d'autres références.

le tome second des « Forces productives et commerciales de la France »[18].

> Remarquez, à partir de Genève jusqu'à Saint-Malo, une ligne tranchée et noirâtre qui sépare le nord et le midi de la France. Au nord, se trouvent trente-deux D$^{ts}$., et treize millions d'habitants; au sud, cinquante-quatre D$^{ts}$., et dix-huit millions d'habitants. Les treize millions d'habitants du nord envoient à l'école 740,846 jeunes gens; les dix-huit millions d'habitants du midi envoient à l'école 375,931 élèves. Il en résulte que sur un million d'habitants, le nord de la France envoie 66,988 enfants à l'école, et le midi, 20,885. Ainsi l'instruction primaire est trois fois plus étendue dans le nord que dans le midi.

Le constat sur les différences dans le nombre d'élèves scolarisés n'était pas nouveau. Il avait été fait par Conrad Malte-Brun (1775–1826), dans son compte-rendu de la statistique européenne de Balbi, publié dans le Journal des Débats du 21 juillet 1823[19].

> Ceux qui croient encore à la possibilité d'administrer les vastes provinces de la France par des institutions uniformes, y verront, par la comparaison du nombre d'élèves qui suivent les académies, les collèges et les écoles dans les diverses parties du royaume, que le midi de la France offre avec le nord le contraste le plus surprenant. Nous savions que la vivacité gasconne et provençale ne s'accommodoit[20] pas d'études méthodiques. Nous savions que l'idiome particulier des Bas-Bretons, la pauvreté de quelques provinces centrales et les mœurs singulières de la Corse, s'opposoient à l'action des lois générales, relatives à l'instruction publique ; mais nous n'aurions jamais imaginé que plusieurs portions de la France méridionale et occidentale fussent placées au niveau des pays européens les plus dépourvus des moyens d'instruction, tandis que sous ce même rapport, le nord de la France marche l'égal des contrées les plus civilisées du monde. Les preuves sont positives ; elles se trouvent réunies dans l'ouvrage de M. Balbi, tome 2, page 246. […] L'instruction publique dans le midi et l'ouest de la France est donc à celle qui existe dans le nord et l'est comme 1 est à 2 1/2.

Alors que l'article de Malte-Brun n'avait pas suscité de réaction particulière, qu'est-ce qui a donné à la conférence de Dupin un tel retentissement ? La visualisation certainement, surtout la symbolique de la couleur « noirâtre » associée à l'obscurantisme ; mais aussi la provocation, dont il n'est pas toujours évident qu'elle soit volontaire. Les « Forces productives de la France », qui ne décrivent en fait que les 32 départements de la France du nord, s'ouvrent sur un « Hommage aux habitants de la France méridionale », qui tient plus du chiffon rouge que de l'hommage, et dont Dupin s'étonne[21] qu'il ait été censuré.

> Je présente à votre émulation généreuse, à votre imitation raisonnée, le modèle d'une partie du royaume favorisée par une longue suite d'événements, favorisée surtout par le voisinage de peuples très-avancés en industrie, et très-heureux en institutions, comme les peuples britanniques, helvétiques et bataves ; tandis que vous n'avez pour voisins que ces peuples d'Espagne et de Portugal, de Sardaigne et d'Afrique, depuis long-temps retardés, dégradés par de mauvaises lois et de mauvais gouvernements.

> Vous serez frappés d'étonnement, lorsque vous verrez quelles différences de population, de richesse territoriale, manufacturière et commerciale, présentent les deux grandes divisions de la France, que nos ancêtres distinguaient en pays de la langue d'oïl et pays de la langue d'oc.[…]

> Déjà quelques-uns de ces résultats vous sont connus ; ils ont causé parmi vous une extrême surprise, et je dirais presque, ils ont excité l'indignation de vos cœurs généreux et fiers.[…]

> Habitants du midi, vous allez connaître, autant que les moyens d'information dont j'ai pu disposer me l'ont permis, l'espace que vous avez à parcourir pour atteindre vos compatriotes du nord. Ils vous ont tracé la route, et l'ont aplanie pour vous[22].

---

18 Dupin, C., 1827, tome II, p. 251 et 337.
19 Essai statistique sur le Royaume de Portugal et d'Algarve, comparé aux autres Etats de l'Europe, par M. Adrien Balbi. III$^e$ article, par Malte-Brun. Adrien-Michel Guerry, qui a collaboré avec Balbi sur le même sujet, ne se prive pas de mentionner l'antériorité de ce dernier en citant l'article de Malte-Brun, dans Guerry, 1833, p. 45.
20 Dans toutes les citations, l'orthographe de l'original a été reproduite.
21 Dupin, C., 1827, tome I, avant propos, p. iv.

On imagine aisément l'indignation qu'un tel discours a pu provoquer parmi les destinataires. Il convient pourtant de ne pas en exagérer la violence. Quelques mois plus tard, le 17 novembre 1827, Dupin était élu député par un des départements noircis sur sa carte, le Tarn. Le 25 décembre, exposant au Conservatoire le même sujet que l'année précédente, Dupin revient sur les réactions que sa carte a suscitées.

> Au premier instant où j'osai mettre en parallèle et rendre sensible à la vue les inégalités frappantes d'industrie et d'instruction populaire, dans le nord et dans le midi de la France, quelques départemens, plus pénétrés que les autres de leur excellence particulière, surtout près d'un fleuve [le Rhône] dont les eaux semblent faire naître ce sentiment, élevèrent le cri de la fierté blessée. On eut dit qu'ils s'en prenaient à moi de n'avoir pas sur ma carte une teinte assez flatteuse, et que le midi tout entier devait à jamais me vouer anathème. Mais au moment où la sagesse du monarque demanda, pour notre pays, des mandataires véridiques avant tout, pour bien servir la patrie par leur bonne foi, au fond du midi, au pied des Cévennes, une ville généreuse [Castres] méprisa ces prétendus griefs, et tournant ses regards vers cet amphithéâtre d'où ne sortirent jamais que des paroles de vérité, la ville industrieuse laissa tomber son suffrage sur celui qui, dans cette enceinte, avait obtenu le vôtre[23].

Plutôt que de reproduire des litanies de plaintes outragées, esquissons une classification des argumentaires en distinguant deux types, anecdotique et méthodologique. Le premier, le plus fréquent, est bien connu des statisticiens qui y sont souvent confrontés. Devant un résultat statistique, le contradicteur multiplie les exemples, qui n'infirment pas pour autant le résultat moyen. Ainsi nombre d'opposants font appel à leurs célébrités locales ; mais Dupin avait par avance contré cet argument en citant lui-même Descartes, Mirabeau, Fermat, et une bonne dizaine d'autres gloires méridionales[24]. Tout aussi faux et fréquent est l'argument consistant à mettre en exergue telle ou telle réalisation. Par exemple le 1er mai 1828, le préfet de la Haute-Loire, à l'occasion de l'« installation des écoles industrielles » au Puy, commence son discours par :

> Lorsqu'un homme d'un vaste savoir, s'appuyant avec confiance sur des renseignemens fournis sans doute avec trop de précipitation ou de légèreté, a cru devoir jeter sur nos contrées le voile de l'obscurité la plus grande, certes il ignorait et le bien que vous avez fait et les succès que vous avez obtenus[25].

Plus rares sont les contradicteurs qui attaquent Dupin sur sa méthode (recueil ou interprétation des données). Dans le même numéro des Annales de la Société d'agriculture, sciences, arts et commerce du Puy, paraît une « Notice sur l'instruction publique, l'Agriculture et l'Industrie de l'arrondissement d'Yssingeaux », par M. de Sainte Colombe. De manière significative, elle est placée sous la rubrique « Statistique » dans la table des matières.

> Le savant auteur des *Forces productives de la France*, en évaluant l'instruction et l'industrie de chaque département, a placé celui de la Haute-Loire dans les degrés inférieurs de l'échelle. La ceinture rembrunie dont il l'a enlacé ressemble à cette lisière qui accuse la faiblesse des premiers pas de l'enfance. Le Velay ne serait-il donc pas sorti des langes du 13e siècle ? Enfans déshérités de cet esprit national sur qui l'émulation exerce un si noble empire ; ses habitans, placés au centre de la France, n'ont-ils rien vu, rien compris, rien essayé de ce qui se faisait autour d'eux? La sombre bannière sous laquelle M. Dupin les a enrôlés pourrait le faire croire ; il n'en est point ainsi cependant ; notre honneur est intéressé à le prouver. […].

---

22 Dupin, C., 1827, dédicace.
23 Dupin, C., 1828, p. 330.
24 Dupin, C., 1827, dédicace, p. v.
25 Annales de la Société d'agriculture, sciences, arts et commerce du Puy, pour 1828. Le Puy, Pasquet, 1829, p. 13.

> C'est par des faits avérés, pouvant défier l'investigation, que nous réclamerons, avec espoir de succès, contre l'arrêt de M. Dupin. C'est en empruntant sa méthode claire, et c'est avec des chiffres que j'établirai mon travail. Les résultats qu'ils présenteront sont plus faciles à saisir et plus propres à amener la conviction. Ils conduiront, j'aime à le croire, au but que je me suis proposé, la réhabilitation de notre pays[26].

Mais la « méthode claire » de Dupin, n'est pas approuvée par tout le monde. La Revue de l'Ouest se permet d'ironiser après coup, sans toutefois citer l'intéressé.

> L'on n'a pas oublié que M. Appert fit, il y a environ un an un voyage en Bretagne : mais, loin d'imiter l'exemple de l'illustre statisticien qui le suivit, ce ne fut ni dans des banquets, ni dans des séances de sociétés académiques que M. Appert recueillit les renseignements qu'il était venu chercher dans le pays[27].

Évoquons maintenant un point de vue moins passionnel, celui des journaux nationaux. La conférence du 29 novembre 1826 est relatée en détail dans les deux journaux à plus fort tirage de l'époque : le Constitutionnel[28] et le Journal des Débats[29]. De larges extraits sont proposés. Le ton est déférent, voire même obséquieux.

> Rien sans doute de plus concluant, rien de plus fort et de plus serré que ce tissu de raisonnemens, tous appuyés sur des faits précis, et dont l'exactitude n'est pas susceptible de la plus légère contestation[30].

Quelques mois plus tard, les « Forces productives et commerciales de la France reçoivent lors de leur parution un large écho de la part des deux journaux. Le Constitutionnel publie deux longs articles[31], le Journal des Débats pas moins de quatre[32]. Si le ton reste globalement favorable, quelques critiques pointent dans le dernier article du Journal des Débats, ainsi qu'une certaine empathie pour la frustration des méridionaux.

> En général, il y a toujours un peu à se défier de la tendance qu'ont les statisticiens pour les contrastes heurtés, le fracas des *chiffres*, la progression gigantesque des nombres, les résultats surprenans. Pourquoi, par exemple; avoir ici établi la division du nord de du midi de la France, par une diagonale, au lieu de la ligne horizontale usitée en géographie pour marquer les latitudes? […]
>
> Les choses ainsi rétablies dans l'ordre qu'il a plu au créateur de leur assigner, et sans rien changer d'ailleurs aux oracles de la statistique, que l'on refasse le compte des *chiffres* de la carte de M. Dupin, et l'on trouvera l'étendue moyenne de l'instruction primaire dans les départemens du midi si peu différente de celle des départemens du nord, qu'il n'y aura plus pour personne sujet de se fâcher. […]
>
> Même en supposant exacts et complets les documens sur lesquels M. Dupin a dressé sa fameuse carte, ma raison se refuse à admettre, sans beaucoup de restrictions, les conséquences qu'il en déduit en termes si absolus ; à mon sens, cette pièce, dont on a fait si grand bruit, n'est pas ce qu'il y a de plus satisfaisant dans le livre auquel elle sert de frontispice[33].

Dupin ne manque pas de répondre longuement à ces critiques[34], augmentant encore la publicité autour de sa carte. Le Figaro, qui est encore à cette époque un journal satirique, ironise sur les réactions de colère, sans la moindre critique à l'égard de Dupin.

---

26 Ibid., p. 80-151.
27 La Revue de l'Ouest, 22 juillet 1829, p. 2.
28 Le Constitutionnel, 11 décembre 1826, p. 2.
29 Journal des Débats, 6 décembre 1826, p. 1.
30 Journal des Débats, 6 décembre 1826, p. 2.
31 Le Constitutionnel, 22 et 26 juillet 1827.
32 Journal des Débats, 24 mai, 9 juin, 1er août, 1er octobre 1827.
33 Journal des Débats, 1er octobre 1827, p. 2-4.
34 Journal des Débats, 4 octobre 1827, p. 2-3.

> Eh ! monsieur le baron Charles Dupin, que vous a donc fait cette pauvre nation française pour que vous la poursuiviez de vos teintes brunes et noires, qui nous prouvent que la moitié de *notre beau pays*, comme le disent les vaudevillistes et les poètes soi-disant patriotiques, est encore plongée dans une ignorance crasse ? On se dit ces choses-là à soi-même ; mais apprendre aux étrangers, par des tableaux fidèles, que les sots et les ignorans sont chez nous en majorité, voilà ce qui est impardonnable. Vos démonstrations sont claires, on ne peut mettre en doute l'évidence de vos preuves. Raison de plus pour armer la colère des départemens noircis sur votre carte accusatrice[35].

Ainsi, la polémique provoquée en France par la carte de Dupin, va donner à l'auteur une position d'autorité incontestée : pour la première fois de son existence, la statistique en tant que discipline est incarnée aux yeux du grand public, et sa popularité s'en trouve décuplée.

### *La statistique morale d'Andraud*

Les érudits de province, qui se sentent tenus de se répondre à Dupin, vont profiter du surcroît de notoriété que la discipline vient d'acquérir grâce à lui. En 1828 et 1829, le mouvement se cristallise autour de l'initiative d'Antoine Andraud (1797–1859), qui est né sur la « ligne noirâtre », à Moulins dans l'Allier. L'idée d'Andraud n'a apparemment rien de polémique :

> La France manque d'une statistique de ses productions intellectuelles ; l'ouvrage que nous offrons au public sous le nom de *Biographie par départemens*, est destiné à remplir cette lacune[36].

Une biographie donc, comme il en fleurit de nombreuses à l'époque ; mais son objectif est ouvertement anti-parisien.

> Jusqu'ici un système funeste de centralisation a pesé sur les provinces. La France n'était que dans Paris, et les autres villes, avec la soumission que donne l'habitude, s'étaient résignées à payer tribut à la capitale qui leur envoyait en échange ses modes, ses livres et ses opinions. A mesure que nous avançons dans le régime constitutionnel, Paris perd insensiblement de ses privilèges, le cercle étroit de la civilisation s'agrandit pour se rompre, chaque ville devient un foyer de lumières indigènes ; des rapports mutuels d'industrie et de savoir s'établissent partout, et cette grande chaîne qui lie nos cités entr'elles ne commence nulle part[37].

De l'aveu même de Dupin, c'est des Bouches-du-Rhône que lui viennent les réactions les plus véhémentes. Il s'y crée, en mars 1827 soit deux ans avant la « Société française de Statistique universelle[38] », une « Société royale de Statistique » dont Dupin sera bientôt membre. Un « Athénée de Marseille » naît en 1828, et dans sa notice historique, il est précisé que « Marseille se trouva ainsi vengée du reproche d'obscurantisme qu'une prévention injuste lui avait adressé[39] ». Tout naturellement, les Bouches-du-Rhône fournissent la première livraison de la biographie par départemens. Lorsqu'elle paraît, et bien que le nom de Dupin n'y soit pas cité, les journalistes comprennent l'objectif.

> Voici venir un correctif utile aux cartes ombrées de M. le baron Charles Dupin. Des noms historiques persuadent mieux que des chiffres au moins incertains, et la reproduction de faits connus inspire un peu plus de confiance qu'une teinte plus ou moins foncée. […]
>
> Comme l'honorable et célèbre professeur-député, M. Andraud (de l'Allier) a entrepris son ouvrage dans le but d'établir parmi nos départemens une louable émulation ; mais nous croyons que l'auteur de la *Statistique morale* a été beaucoup mieux inspiré que l'orateur du Conservatoire des arts et métier. Ce dernier n'a rien trouvé de mieux que de noircir les départemens qui lui semblaient être le plus en arrière dans la route de la civilisation ; M. Andraud veut leur inspirer le goût des sciences et des arts, en plaçant sous leurs yeux les titres et la gloire de leurs compatriotes. M. le baron Dupin s'est contenté de blesser l'amour-propre des hommes, tandis que son compétiteur excite en eux un noble orgueil par de beaux exemples[40].

---
35  Le Figaro, 10 juillet 1828, p. 686.
36  Andraud, A., 1829, préface.
37  Ibid.
38  De Falguerolles, A. 2010.
39  Comte, A.-J., 1846, p. 406.

L'idée lumineuse d'Andraud, est d'avoir fait des souscripteurs les rédacteurs de son ouvrage. L'appel à souscription, largement diffusé dans la presse nationale en 1828, dévoile l'astuce :

> *Nota* : Les personnes qui ont des droits quelconques à la reconnaissance et à l'estime de leurs concitoyens, et qui désirent voir figurer leurs noms dans la Biographie de leur département, sont priées d'adresser leurs notes et renseignemens *franco* au bureau de la direction, rue de Richelieu, nº 67, à Paris[41].

Voir son nom imprimé parmi ceux des gloires locales ! Qui aurait résisté à l'achat de quelques volumes à offrir aux amis et connaissances ? Le succès commercial était assuré, quitte à ce que la sélection ne soit pas des plus rigoureuses. Quelques pamphlétaires ne se privent pas de caricaturer ces « nullités avides de renommée[42] », en citant des noms.

> Le lecteur même qui connaît le département ne sera-t-il pas stupéfait quand il verra trois colonnes de ce recueil destinées à apologiser un annotateur, Julliany, dont je n'aurais jamais soupçonné l'éclatant mérite sans les phrases ridicules d'éloges qu'un autre lui jeta, dit-on, au visage dans une certaine Société statistique aussi inconnue que ses honorables membres[43] ?

Mais globalement la presse, lors de l'appel à souscription, puis de la parution des premiers numéros à laquelle elle donne un large écho, reste très favorable à l'initiative.

> Quatre livraisons de cette importante Collection ont déjà paru ; elles suffisent pour faire juger l'ouvrage entier, et garantissent le succès d'une entreprise qui doit puiser, dans la faveur même avec laquelle le public l'accueille, de nouveaux élémens de perfection. Les Bouches-du-Rhône, le Var, les Basses-Alpes et le Gard, placés maintenant sous les yeux des lecteurs, déroulent le brillant tableau des illustrations méridionales. […]
>
> Les livraisons de l'*Hérault* et de *Vaucluse*, dont la librairie annonce la prochaine publication, complèteront le tableau biographique des départemens du Midi. On ne peut qu'encourager l'éditeur à persévérer dans l'exécution de la patriotique pensée ; le public, qui n'a pas attendu ces éloges pour lui accorder sa faveur, ne se lassera pas de prêter son appui et d'applaudir à cette grande entreprise[44].

Pour une raison inconnue, la « grande entreprise » n'a pas dépassé les quatre premières livraisons. Le biographe d'Andraud estime que « la révolution de 1830 arrêta à son début cette grande publication[45] », mais cette explication paraît peu compatible avec le rythme rapide des premières parutions, au printemps 1829.

Même si l'on reste conscient des biais induits par le mode de sélection, les premiers volumes de la « Statistique Morale » d'Andraud fournissent un échantillon instructif des « illustrations méridionales ». Il est intéressant de le comparer, d'une part avec la répartition des collaborateurs à la statistique des départements sous le consulat établie par Perrot[46], d'autre part avec les professions représentées dans une société savante de province, au Puy, pour la même année 1829[47]. Ce qui frappe parmi les 296 noms de la biographie d'Andraud pour les Bouches-du-Rhône, c'est la sur-

---

représentation des professions juridiques (juges, avocats, notaires, etc.), qui constituent 21 % du total, contre 10 % dans la société savante du Puy. Perrot compte 4 % de magistrats, ecclésiastiques, et officiers, tandis que Andraud liste déjà 5 % d'ecclésiastiques et 9 % de militaires. En revanche, les professions médicales ne constituent que 9 % du total chez Andraud, mais 14 % chez Perrot, et 16 % au Puy. Une catégorie est complètement absente dans la biographie d'Andraud, celle des propriétaires, qui sont 23 % chez Perrot, 27 % au Puy. Cette particularité s'explique par l'examen des notices. Nombre de propriétaires, auteurs de diverses pièces écrites, préfèrent se présenter comme « littérateurs » : ils sont 15 % chez Andraud. Autre particularité, 12 % des notices d'Andraud relèvent de professions artistiques (peinture, dessin, musique, théâtre, etc.). Il est compréhensible que les artistes n'apparaissent pas dans le décompte de Perrot, et peu au Puy (4 %). La proportion des négociants ou chefs d'entreprise, qui n'était que de 4 % chez Perrot, passe à 10 % chez Andraud, et 9 % au Puy. Globalement, la variabilité observée s'explique par les différences de définition des populations : estimer avoir droit à « la reconnaissance et l'estime de ses concitoyens », ne signifie pas que l'on doive produire des travaux d'ordre statistique, ni être membre d'une société savante. Une constante pourtant, ne surprendra personne : la quasi absence de femmes. Seulement trois sont citées pour leurs mérites personnels par Andraud. L'évolution des mentalités était pourtant déjà engagée, comme en témoigne le compte-rendu de la séance du 20 décembre 1836 à la Société française de Statistique universelle[48]. Une proposition, déposée par le comte Legrand, était ainsi rédigée.

    1° Les dames ne seront plus reçues à l'avenir dans la Société française de Statistique universelle;

    2° Les cinq dames reçues jusqu'à ce jour, continueront cependant à en faire partie.

Reconnaissons à nos prédécesseurs l'honneur d'avoir, après de longs débats passionnés, rejeté la proposition à la quasi unanimité[49].

## 2. Un développement exubérant

*Statistique journalistique*

Vu le sujet de cet article, il aurait été pour le moins paradoxal qu'une tentative de validation statistique n'y soit pas proposée ; la voici. Les 12612 numéros du Journal des Débats parus entre 1814 et 1848, ont été choisis comme corpus de référence[50]. Pour chacune des 35 années de parution, la proportion des numéros de l'année contenant les chaînes de caractères « statistique », « Charles Dupin » ou les deux, a été relevée. L'étude a été complétée par un examen individuel des occurrences. Les résultats numériques sont représentés sur la figure 1. En comparant les périodes 1814-1826 et 1828-1848, on y constate une nette augmentation de la moyenne des trois proportions, que le test de Student confirme (p-valeurs $1,81 \times 10^{-21}$, $1,02 \times 10^{-11}$, $1,71 \times 10^{-8}$). Ce résultat statistique brut doit être interprété avec beaucoup de prudence, et les choix méthodologiques qui y conduisent sont tous discutables. Nous examinerons successivement le journal, les chaînes de caractères, et les indicateurs.

Concernant le journal, la restriction technique de la disponibilité en mode texte est une contrainte importante. Elle écarte entre autres le Constitutionnel, qui a le plus fort tirage de la presse de l'époque. Le Figaro est disponible en mode texte, mais les résultats y sont moins significatifs. En premier lieu, le Figaro est de création récente (1826), et la comparaison avec la période antérieure est impossible. D'autre part, il s'agit à ses débuts d'un journal satirique dans lequel on trouve moins d'articles généralistes, et surtout moins d'annonces de parution que dans le Journal des Débats. Les résultats concernant Charles Dupin y font apparaître une anomalie : du printemps à l'automne 1831, celui-ci a été la cible d'une intense campagne de dénigrement de la part du Figaro. Ses prestations oratoires, ses déboires électoraux, ses opinions politiques variables, son goût pour le cumul de fonctions, de rémunérations et de titres, et bien sûr son attrait pour les statistiques et les cartes ombrées, y sont très souvent caricaturés. L'extrait suivant donne une idée du ton.

---

48 Journal des travaux de la Société française de Statistique universelle, 2(19), janvier 1837, p. 386.
49 On trouvera dans les noms des premières statisticiennes dans de Falguerolle, A., 2010, p. 17.
50 Accessible en mode texte sur gallica.bnf.fr.

> Ensuite on a trouvé défunt tout pareillement M. Charles Dupin, baron par quittance, bariolé de mathématiques, orateur fatigant et infatigable, rongé d'impopularité, rêvasseur de fictions statistiques, et badigeonné d'encre de Chine, l'homme de l'ubiquité et du cumul[51].

Le Journal des Débats, s'il offre plus de garanties d'objectivité que le Figaro vis-à-vis de la statistique et de Charles Dupin, n'est pas pour autant exempt de tout biais politique. D'abord très dévoué à la Restauration[52], le journal, tout en restant monarchiste et catholique, évolue vers une opposition modérée sous le ministère Villèle (1821–1828), auquel il « déclare la guerre » à partir de 1826[53]. Il n'est pas possible d'exclure d'emblée que l'augmentation de fréquence relevée à partir de 1828 puisse être liée au changement de ligne politique du journal.

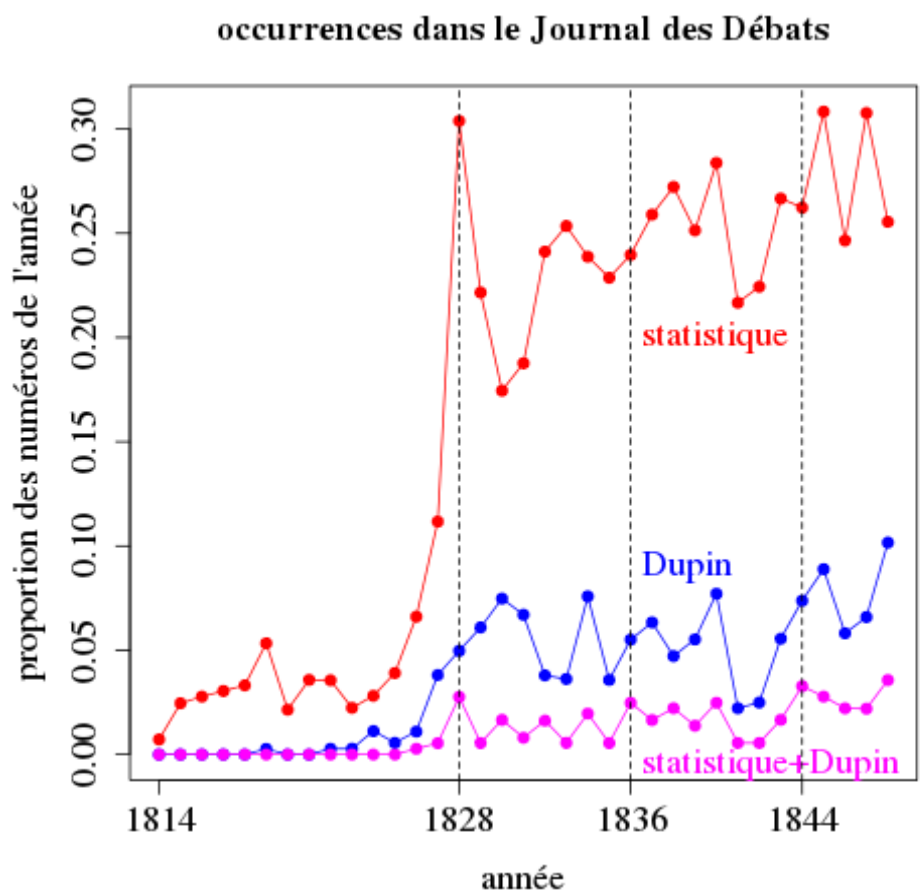

*Figure 1 : proportions annuelles des chaînes de caractères « statistique », « Charles Dupin » ou les deux dans le Journal des Débats, de 1814 à 1848. Les années où la troisième est significativement élevée par rapport aux deux premières sont marquées par des lignes verticales.*

---

51 Le Figaro, 14 juillet 1831, p. 2.
52 Nettement, A., 1838, p. 13.
53 Nettement, A., 1838, p. 57.

Concernant les chaînes de caractères, on doit d'abord observer que la recherche automatique induit un biais négatif : certaines occurrences peuvent ne pas avoir été détectées, et ce d'autant plus que la chaîne cherchée est plus longue. Plutôt que « statistique », on aurait pu rechercher « statisti », qui inclut également « statisticien » et « statistiquer ». Le nom statisticien est utilisé dès 1805 par Donnant[54], on le trouve dans la presse en 1810[55], mais il reste peu fréquent : seulement 43 numéros du Journal des Débats l'utilisent entre 1828 et 1848 contre 1500 pour statistique. Le verbe statistiquer, employé par Balzac en 1830, est encore plus rare. Utilisé dans une revue de vulgarisation en 1862[56], on ne le trouve dans la presse généraliste qu'après 1880[57]. Cela peut-être lié au succès de l'opérette « la belle Lurette » et ses « Couplets de la statistique[58] ». La recherche concernant Dupin pose un problème d'homonymie. On trouve selon les années, entre 6 et 7 fois plus d'occurrences pour « Dupin » que pour « Charles Dupin », en grande partie à cause des deux frères de Charles : André (1783–1865), et Philippe (1795–1846). Un autre baron Dupin, Claude François Étienne (1767–1828), lauréat du prix Montyon de statistique en 1822, apparaît aussi sur la période[59] ; il n'est qu'un lointain parent des précédents. Rechercher « Charles Dupin » résout le problème d'homonymie mais aggrave la sous-estimation, car il est souvent désigné par « M. le baron Dupin », ou bien « M. Ch. Dupin ». L'examen des occurrences du mot statistique donne une bonne idée de l'extension foisonnante du thème : voici une liste alphabétique d'adjectifs qui lui sont accolés comme épithètes.

> agricole, agronomique, archéologique, botanique, bourgeoise, bovine, carcérale, chimique, chinoise, civile, commerciale, conjugale[60], criminelle, ecclésiastique, électorale, équestre, européenne, forestière, géologique, historique, industrielle, judiciaire, littéraire, matrimoniale, médicale, militaire, minéralogique, morale, mortuaire, monumentale, musicale, nationale, nosologique, œnologique, parisienne, philosophique, picturale, pittoresque, politique, rurale, sociale, théâtrale, topographique, universelle, végétale, zoologique.

Le choix de la proportion annuelle des numéros comme indicateur, est également discutable. On aurait pu baser l'étude sur le nombre total d'occurrences, ou la distribution des nombres d'occurrences par numéro, mais l'indicateur choisi est moins sensible à la sous-estimation par l'algorithme de reconnaissance. Calculer une proportion plutôt que de conserver les nombres bruts est motivé par le fait que le nombre annuel de numéros est plus variable qu'il n'y paraît : de 277 en 1814 à 375 en 1819, avec les numéros spéciaux. Il se trouve que dans la grande majorité des cas, le mot statistique s'il est rencontré dans un numéro, n'y figure qu'une fois ; le plus souvent, il s'agit de l'annonce de la parution d'un ouvrage, dont le titre ou la description contient le mot. Par exemple sur les 110 numéros comportant le mot statistique en 1828, celui-ci n'apparaît que dans 17 articles de fond, et seuls 15 numéros contiennent le mot plus d'une fois. Puisque les occurrences sont très majoritairement des citations, il semble que le biais politique évoqué plus haut puisse être écarté : l'augmentation observée est essentiellement due aux nouvelles parutions. Par exemple, toutes les occurrences isolées de janvier à mars 1828 se trouvent dans les annonces successives de « cartes séparées de chaque département », intitulées « Tableau géographique et statistique ». C'est bien parce que le mot statistique est associé à des descriptions, ou des guides, géographiques et/ou historiques le plus souvent, qu'il devient plus fréquent. Par exemple, en avril 1828 un « Manuel général et nouveau des voyageurs, négocians, administrateurs, hommes d'affaires », est sous-titré

---

54  Donnant, D.-F., 1805, utilise le mot statisticien à 21 reprises.
55  Magasin encyclopédique ou journal des sciences des lettres et des arts, tome III, 1810, p. 430.
56  La Science pittoresque, 26 décembre 1862, p. 418.
57  La Presse, 18 mars 1880, le Gaulois, 9 avril 1882, le Figaro, 7 novembre 1883, Gil Blas, 8 mai 1889, etc.
58  Offenbach, J., et al., 1880, p. 38.
59  Journal des Débats, 15 avril 1822, p. 4.
60  Chez Balzac, H. de, 1846, à qui l'on pardonnera cette légère incohérence.

« Statistique complète de la France[61] ».

Reste à examiner la relation avec les références à Charles Dupin. La figure 1 semble indiquer que les occurrences de « statistique » et « Charles Dupin », sont très corrélées et c'est effectivement le cas (coefficient : 0,86, p-valeur : $1,43 \times 10^{-11}$). Mais que les occurrences annuelles soient corrélées ne signifie pas qu'elles soient souvent simultanées. Si on teste la taille des occurrences simultanées, on trouve qu'elle n'est significative que pour trois années, 1828, 1836 et 1844 (p-valeurs :  0,020, 0,028, 0,026). Or il se trouve que la plupart de ces occurrences communes sont en fait séparées : les deux termes apparaissent dans des articles différents du même numéro. Par exemple en 1828, sur les 10 numéros qui contiennent à la fois statistique et Dupin, un seul les associe effectivement[62]. Il semble donc légitime de conclure que, malgré la stature de spécialiste que Dupin a acquise auprès du grand public, l'augmentation de fréquence du mot statistique n'est pas liée à l'accroissement de sa propre notoriété.

*Personnages de théâtre et de roman*

La scène se passe dans un jardin. Édouard et Amélie espèrent se marier. Caché derrière un arbre, le marquis de Morville, oncle d'Édouard, épie les amoureux.

> AMÉLIE : Hélas ! quand toutes ces divisions finiront-elles ?
>
> ÉDOUARD : Bientôt, heureusement pour la France. Vous connaissez sans doute les nouveaux calculs statistiques ?
>
> LE MARQUIS : De la statistique ! Oh ! le nigaud !
>
> ÉDOUARD : Ils sont très rassurants : M. Dupin prouve fort bien…
>
> LE MARQUIS : Qu'il est un cannibale.
>
> ÉDOUARD : Que, dans quelques années, toute la vieille génération aura disparu. Nous la remplacerons au pouvoir.

La pièce, « les Stationnaires », est datée de septembre 1827. Elle est la dernière du second volume des « Soirées de Neuilly », par M. de Fongeray[63]. Les comptes-rendus dans la presse sont élogieux, et soulignent le personnage d'Édouard « bien comique et bien vrai, quoiqu'il ne semble pas d'abord le plus gai de la pièce[64] ». Stendhal précise :

> la scène la plus comique de la pièce est celle qui ridiculise le sérieux affecté à la mode aujourd'hui parmi nos jeunes gens. L'auteur fait le portrait d'un grave jeune homme qui, au lieu de faire la cour à la jeune fille qui doit être sa femme, lui fait un long discours sur les statistiques. Il fallait un certain courage pour se moquer ainsi d'un travers aussi répandu, car les jeunes gens que l'auteur a ridiculisés collaborent à nos principaux journaux littéraires, et Fongeray y sera sans doute traité d'écrivain excessivement immoral[65].

Dès septembre 1827, le personnage du statisticien trop sérieux, qui débite des évidences (« toute la vieille génération aura disparu ») en citant Dupin (qui est un cannibale), est donc suffisamment reconnu par le public pour être utilisé comme ressort comique. Retrouvons-le trente-quatre ans plus tard, chez Eugène Labiche (1815–1888). Il s'appelle Magis et courtise Lucile qui est la fille de Madame de Guy. Il est présenté par Désambois ; Horace est le témoin sceptique.

---

61  Journal des Débats, 25 avril 1828, p. 4.
62  Journal des Débats, 24 juin 1828, p. 2.
63  De Fongeray, J. F., 1828, p. 415. De Fongeray est un pseudonyme. Les auteurs sont Hygin-Auguste Cavé (1796–1852) et Adolphe Dittmer (1795–1846) : voir  Asse, E., 1900, p. 271.
64  Journal des Débats, 11 août 1828, p. 3.
65  Stendhal, 1936, p. 331 . La lettre est datée du 21 mars 1828. Stendhal avait posé pour le « portrait de l'auteur » des Soirées de Neuilly : Asse, E., 1900, p. 272.

DÉSAMBOIS, *saluant* : Monsieur. *(Aux dames.)* Je précède de quelques instants M. Célestin Magis, l'heureux compétiteur à la main de mademoiselle. Croyez bien que je ne lui aurais pas accordé mon patronage, si je n'avais distingué en lui les qualités les plus solides… M. Magis est un jeune homme sérieux… tout à fait sérieux !

HORACE, *à part* : Encore ? Ah çà ! ils sont donc tous sérieux ?

MADAME DE GUY : Vous savez, monsieur Désambois, que j'ai toute confiance en vous !

HORACE : Peut-on, sans indiscrétion, demander quelle est la profession de mon futur cousin ?

DÉSAMBOIS : Mon Dieu, il n'en a pas positivement… c'est un homme…

HORACE : Sérieux ?…

DÉSAMBOIS : Oui… qui s'occupe de sciences… d'études transcendantes !

HORACE : Ah !

DÉSAMBOIS : A vingt-neuf ans, M. Magis vient d'être promu à la dignité de secrétaire de la Société de statistique… de Vierzon…

HORACE : Fichtre ! c'est un beau grade !

DÉSAMBOIS : Et j'ai moi-même l'espoir d'être avant peu nommé membre correspondant de ladite…

[…]

DÉSAMBOIS : Si j'avais un fils, je voudrais, qu'il vous ressemblât… Monsieur Magis, il faudra envoyer à ces dames votre dernier ouvrage. *(A Lucile.)* Il a publié un ouvrage… imprimé…

MADAME DE GUY : Comment ?

MAGIS : Je n'aurais pas osé prendre cette liberté ; mais, puisque vous le permettez, je serai heureux de vous apporter moi-même mon opuscule sur la *Monographie de la statistique comparée*.

DÉSAMBOIS : Avec un petit mot sur la première page…

MADAME DE GUY : Ah ! monsieur, vous ne pouvez douter de l'intérêt…

DÉSAMBOIS, *bas* : Laissez-le parler !

MAGIS : La statistique, madame, est une science moderne et positive. Elle met en lumière les faits les plus obscurs. Ainsi, dernièrement, grâce à des recherches laborieuses, nous sommes arrivés à connaître le nombre exact des veuves qui ont passé sur le pont Neuf pendant le cours de l'année 1860[66].

On trouve chez Magis l'ensemble des ressorts comiques du personnage : le statisticien est trop sérieux, ses comptes sont inutiles (le nombre exact des veuves), on ne sait quelle est sa fonction sociale exacte (les études transcendantes), son ambition et sa vanité passent à la fois par la reconnaissance académique (la société de statistique de Vierzon), et l'ouvrage imprimé (la Monographie de la statistique comparée).

Nous avons vu Honoré de Balzac (1799–1850), vitupérer contre les « cacographes » en 1830. En mars 1832, il ironisait sur « certains hommes toujours en travail d'une œuvre inconnue : statisticiens tenus pour profonds sur la foi de calculs qu'ils se gardent bien de publier…[67] ». On s'attendait donc à ce que le personnage figure en bonne place dans la Comédie Humaine. Pourtant Balzac ne fait que l'esquisser, en 1841 dans « Ursule Mirouët ».

> Je dois, dit-il, me faire oublier pendant trois ou quatre ans, et chercher une carrière. Peut-être me ferai-je un nom par un livre de haute politique ou de statistique morale, par quelque traité sur une des grandes questions actuelles[68].

---

66 Labiche, E., 1896, p. 423-431. La première a eu lieu le 16 mars 1861.
67 Balzac, H. de, 1852, p. 293.
68 Balzac, H. de, 1855, p. 102.

Gustave Flaubert (1821–1880) va beaucoup plus loin dans la satire mordante. Son érudit ambitieux est Homais, le pharmacien d'Yonville, qui selon Zola « est l'importance provinciale, la science de canton, la bêtise satisfaite de tout un pays[69] ».

> Cependant, il étouffait dans les limites étroites du journalisme, et bientôt il lui fallut le livre, l'ouvrage ! Alors il composa une *Statistique générale du canton d'Yonville, suivie d'observations climatologiques*, et la statistique le poussa vers la philosophie. […]
>
> Il était le plus heureux des pères, le plus fortuné des hommes. Erreur ! une ambition sourde le rongeait : Homais désirait la croix. Les titres ne lui manquaient point : 1° S'être, lors du choléra, signalé par un dévouement sans bornes ; 2° avoir publié, et à mes frais, différents ouvrages d'utilité publique, tels que… (et il rappelait son mémoire intitulé : *Du cidre, de sa fabrication et de ses effets* ; plus des observations sur le puceron laniger, envoyées à l'Académie ; son volume de statistique, et jusqu'à sa thèse de pharmacien) ; sans compter que je suis membre de plusieurs sociétés savantes (il l'était d'une seule)[70].

Émile Zola (1840–1902) nous a laissé ce portrait savoureux de six « membres distingués de la société de statistique » : « grands pieds, grandes mains, larges figures massives », qui, « muets, approuvaient tout de la tête[71] ». Pourtant, plutôt qu'un ambitieux, l'amateur de statistique est chez Zola un compilateur obsessionnel, qui peut pousser sa manie jusqu'à la folie ; comme Mouret qui « compte les s qui se trouvent dans la Bible[72] », Vabre qui depuis dix ans, « dépouillait chaque année le catalogue officiel du Salon de peinture, portant sur des fiches, à chaque nom de peintre, les tableaux exposés[73] » ; ou bien cet inconnu rencontré sur le trottoir :

> Ah ! C'est vous, monsieur, me dit-il en balbutiant. Vous devriez bien m'aider à compter les étoiles. J'en ai déjà trouvé plusieurs millions, mais je crains d'en oublier quelqu'une. C'est de la statistique seule, monsieur, que dépend le bonheur de l'humanité[74].

La « Belle Lurette », représentée pour la première fois le 30 octobre 1880, est la dernière opérette de Jacques Offenbach (1819–1880). Malicorne, intendant du duc de Marly, y est dépeint comme un intrigant, mais son amour de la statistique, ridiculement proclamé dans les célèbres « Couplets de la statistique », n'est pas particulièrement lié à son ambition personnelle. L'accent est plutôt mis sur deux autres ressorts comiques : les doutes sur le contour et la définition (« un champ si vaste et si profond que même ceux là qui le font, n'ont jamais pu rien y comprendre »), et la futilité des comptes : (« combien la porte Saint-Denis voit passer de bêtes à cornes »[75]).

*Perceptions de la discipline*

Dès ses débuts napoléoniens, la statistique en France s'est construite autour de diverses oppositions et polémiques, portant parfois sur ses fondements même[76]. Celles-ci sont suffisamment bien connues pour qu'il ne soit pas utile d'y revenir. Nous allons plutôt nous attacher à la perception positive ou négative de la discipline, telle qu'elle se manifeste dans la presse et la littérature. On trouve dans les écrits de l'époque, beaucoup plus d'enthousiasme que de critiques, plus de déférence que de sarcasmes. Avant la période qui nous intéresse ici, il convient tout d'abord de mentionner quelques articles du Journal des Débats, remarquables par leur clairvoyance. Les deux premiers sont présentés comme un exposé de vulgarisation, à l'occasion de la parution d'une « théorie de l'économie politique », par Ch. Ganilh. On trouve dans celui du 18 décembre 1815 un véritable réquisitoire, qui se conclut par :

> retenons bien que la science dite Economie politique est le matérialisme de l'administration ; que cette science se partage en deux sectes, dont l'une procède par imagination, et l'autre par expérience ; mais que l'expérience repose sur la statistique, et que la statistique ne repose sur rien[77].

---

[69] Zola, E., 1906, p. 279.
[70] Flaubert, G., 1879, p. 381-383.
[71] Zola, É, 1876, p. 302-305. Les statisticiens dans les Rougon-Macquart ont été recensés par Brunet, É., 1985.
[72] Zola, É, 1874, p. 293.
[73] Zola, É, 1882, p. 96.
[74] Zola, É, 1864, p. 29.
[75] Offenbach, J., et al., 1880, p. 38. Les librettistes sont Ernest Blum (1836–1905), Édouard Blau (1836–1906) et Raoul Toché (1850–1895).
[76] En particulier le conflit entre Duvillard et Deferrière en 1806 : voir Gille, B, 1980, p. 123.

La sentence sert aussi d'introduction au second article[78] ; elle est suffisamment forte pour avoir été remarquée par Saint-Simon[79]. Vu la période, il serait tentant de mettre l'opinion de l'auteur sur le compte du zèle anti-bonapartiste du journal. Mais il se trouve que cet auteur est Joseph Fiévée[80] (1767–1839), «lequel n'est pas du tout caché sous les lettres T. L.[81] ». Or s'il est effectivement monarchiste et surtout anti-jacobin, Fiévée a eu une expérience de première main, en tant que préfet de la Nièvre du 8 avril 1813 au 22 mars 1815. Il n'a pas publié de statistique de son département, mais en a écrit une description qui aurait pu en tenir lieu[82]. Son argumentaire mérite donc d'être examiné. Il est basé sur la difficulté du recueil de données : pour lui, « toute demande en renseignemens faite par l'autorité inspire toujours des craintes ou des espérances ; et les réponses qu'elle reçoit sont erronées en plus ou en moins, selon que celui qui répond a laissé tourner son imagination vers la crainte ou vers l'espérance ». Il dispose d'exemples précis :

> Dans un moment où l'on craignoit une disette, j'ai vu deux recensemens des grains faits à la même époque dans le même département, l'un par les maires et par ordre du préfet, l'autre par la gendarmerie et sur l'ordre direct de la police. Le recensement fait par les gendarmes offroit un déficit de moitié sur celui fait par les maires ; et je ne crois pas que celui-ci fut exempt d'erreurs volontaires.

Quand bien même n'y aurait-il aucun intérêt en jeu,

> croyez-vous […] que la négligence ne produise pas autant d'erreurs que l'intérêt ? […] D'ailleurs, qui a la certitude que, dans les préfectures, on mette beaucoup d'exactitude à faire un tableau général des tableaux envoyés par les sous-préfets ? Pour moi, je ne répondrois pas que, pour s'éviter des correspondances d'autant plus actives qu'elles n'ont aucun but, et des reproches de négligence inévitables, certains préfets ne fissent dresser leurs tableaux avant les renseignemens, comme cet historien fit pour le siège qu'il avoit à décrire.

La leçon n'avait certainement pas été entendue par Dupin, qui n'émet jamais le moindre doute sur l'« exactitude positive » de ses données. Son assurance est largement véhiculée dans la presse et la littérature. On ne doute pas de la véracité des chiffres : « c'est brutal une statistique, mais c'est net et clair[83] », ou encore « je ne sais rien de plus éloquent que les chiffres[84] ». Peut-être, mais si leur exactitude n'est que rarement remise en cause, leur utilité pour le vivant, en particulier l'humain, est sujette à caution. Par exemple le Tintamarre parle d'une « statistique à la Charles Dupin, c'est-à-dire simple, concise, sèche, absurde[85] ». Auguste Comte (1798–1857), qui récuse toute application des mathématiques à la biologie au nom de la diversité du vivant, n'accorde aucun crédit à la statistique médicale.

---

77  Journal des Débats, 18 décembre 1815, p. 4.
78  Journal des Débats, 9 janvier 1816, p. 3. Fiévée revient sur le sujet le 27 janvier et le 23 mars.
79  Saint-Simon, H., 1817, p. 132.
80  Sur le rôle de Fiévée sous la restauration, voir Popkin, J. D., 2001.
81  Hatin, E., 1861, p. 195.
82  Thuillier, G., 1961.
83  Alexandre Dumas fils (1824–1895) : Dumas, A., 1880, p. 127.
84  Verne, J., 1868, p. 299.
85  Le Tintamarre, 5 janvier 1852, p. 5.

À la vérité, l'esprit de calcul tend de nos jours à s'introduire dans cette étude, surtout en ce qui concerne les questions médicales, par une voie beaucoup moins directe, sous une forme plus spécieuse, et avec des prétentions infiniment plus modestes. Je veux parler principalement de cette prétendue application de ce qu'on appelle la statistique à la médecine, dont plusieurs savans attendent des merveilles, et qui pourtant ne saurait aboutir, par sa nature, qu'à une profonde dégénération directe de l'art médical, dès lors réduit à d'aveugles dénombremens. Une telle méthode, s'il est permis de lui accorder ce nom, ne serait réellement autre chose que l'empirisme absolu, déguisé sous de frivoles apparences mathématiques[86].

En 1822, à l'occasion de la remise du prix Montyon de statistique à un ouvrage intitulé « Observations géognosiques faites dans les Pyrénées », Fiévée s'inquiètait déjà de l'étendue que la discipline était en train de prendre.

On se demande, avec étonnement, depuis quand une science toute entière, la *géognosie*, est devenue partie intégrante de la soi-disant *statistique ?* L'étonnement s'accroît lorsqu'on lit le programme obscur et vague où l'Académie essaie de définir ce terme nouveau de *statistique*. Il y est dit entre autres que « la statistique décrit le climat, le territoire et les divisions politiques ou naturelles ». Voilà la géographie englobée dans la statistique. Elle indique encore « les monumens de l'histoire et des arts » ; la voilà enrichie de tout le domaine de l'archéologie[87].

L'étonnement de Fiévée n'a pas été partagé par les fondateurs de la Société royale de Statistique de Marseille, qui divisent leurs travaux en trois classes, « 1° sciences morales, philosophiques et industrielles ; 2° sciences naturelles, physiques et mathématiques ; 3° sciences qui traitent des langues, de la littérature et des beaux-arts » ; et d'annoncer fièrement quelques années après que « la météorologie, l'instruction publique, la géologie, l'histoire, l'archéologie, la population, les consommations, l'industrie agricole, manufacturière et commerciale, la navigation, etc., ont été étudiées avec soin[88] ». Une telle prétention à l'universalité ne pouvait que susciter l'incrédulité. On ne s'étonnera donc pas que Paul Féval (1816–1887) fasse dire à un des personnages, sur le ton de la provocation, que « plusieurs bons esprits regardent [la statistique] comme devant remplacer toutes les autres [sciences] dans un temps donné[89] ». Il y a eu effectivement de ces « bons esprits », comme le préfet Thiessé, qui s'adressant aux membres de la société de statistique des Deux-Sèvres, leur dit : « Vous vous êtes voués à la statistique, base de toutes les connaissances humaines, introduction nécessaire à toutes les sciences[90] ».

Nombreux sont ceux qui prennent le contre-pied de l'esprit du temps en affectant de ne pas prendre au sérieux cette science trop sérieuse ; comme Edmond About : « c'est moins une science qu'un tâtonnement raisonné[91] », ou Balzac encore qui parle d'« enfantillage des hommes d'Etat modernes[92] ». L'idée que la statistique peut argumenter sur tout et son contraire, se répand : quand un des personnages de Maupassant s'exclame : « Ah ! les statistiques ! On leur fait dire ce que l'on veut, aux statistiques[93] ! », il ne fait que reprendre un jugement maintes fois répété, et qui semble avoir été largement partagé. En 1843, la « revue des deux mondes » publie un long article signé Louis Reyba, intitulé « La société et le socialisme ». La statistique y est accusée de servir la cause des utopies et d'aligner avec « candeur », des calculs qu'elle interprète avec « naïveté » ; mais surtout,

Si la statistique ne sait pas mieux se contenir, elle se fera, auprès des esprits sérieux, un tort irréparable. C'est une science qui renferme des calculs et des argumens pour toutes les causes, fussent-elles diamétralement opposées. Les chiffres sont complaisans ; ils se prêtent aux désirs secrets de l'observateur et à la fortune des livres. On se propose de prouver une chose, et l'on voit tout dans le sens de cette démonstration[94].

---

86  Comte, A., 1838, p. 418.
87  Journal des Débats, 15 avril 1822, p. 4.
88  Comte, A., 1846, p. 420.
89  Féval, P., 1877, p. 57.
90  Revue littéraire de l'ouest, 1837, p. 226.
91  About, E., 1866, p. 143.
92  Balzac, H. de, 1838, p. 331.
93  Maupassant, G. de, 1891, p. 149.

Plus légèrement, Alphonse de Lamartine (1790–1869) est sûr de mettre les rieurs de son côté quand il s'exclame à la tribune de l'Assemblée, le 22 avril 1846 :

> J'ai appris dans cette étude ce que valent les statistiques, j'en demande bien pardon aux partisans de la statistique qui se trouvent ici, excusez cette comparaison vulgaire : les prestidigitateurs font leurs tours avec des gobelets ; les économistes font leurs théories avec des statistiques[95].

Que la statistique produise les décomptes les plus incongrus ou inutiles est souvent raillé, comme nous l'avons vu chez Zola. « La statistique […] fait souvent des merveilles en nous apprenant, par exemple, que la longueur totale des asperges mangées par un homme bien constitué pendant son existence est juste le double de celle du câble transatlantique[96] ». Victor Hugo (1802–1885), n'a pas beaucoup de respect pour ce « faiseur de statistique [qui] a calculé qu'en superposant l'un à l'autre tous les volumes sortis de la presse depuis Gutemberg on comblerait l'intervalle de la Terre à la Lune[97] ». Pour susciter ironie et critiques, les « faiseurs de statistiques » étaient-ils donc si nombreux et enthousiastes ? On serait tenté de croire à une irrésistible lame de fond, à lire le Bulletin de la Société française de Statistique universelle, qui annonce en avril 1830, 441 membres, dont « 5 princes du sang de France et d'Angleterre, 47 pairs du royaume, 24 députés, 12 ministres actuels ou d'état, 20 ambassadeurs ou ministres étrangers, 12 membres de l'Institut, 250 personnes prises dans toutes les notabilités sociales de la France, 71 étrangers[98] ». Mais une lettre de Benoît d'Azy à son ami Félicité de Lammenais (1782–1854), montre l'envers du décor ; elle est datée du 6 mai 1830.

> Je n'ai point encore envoyé ta lettre à M. César Moreau, l'homme aux statistiques, je crois même que tu feras bien de ne pas la lui adresser du tout. Cette prétendue société de statistique s'est formée tout entière dans le cerveau de M. César Moreau homme peu considéré à ce qu'on assure. Il a commencé par envoyer des diplômes à tout le monde et quelques personnes honorables se sont laissé prendre à cette idée d'une chose utile en elle-même ; appuyé sur les premiers noms, il est devenu plus assuré et a prévenu beaucoup de gens que la société les avait nommés *membres ;* tout cela est encore de son imagination car je crois que deux personnes ne se sont encore jamais réunies pour cet objet. J'ai pris quelques renseignements parce que je voyais sur ses feuilles des noms honorables, je n'y ai vu jusqu'à présent que ce que les anglais appellent un *catch penny* et je me suis décidé à attendre. J'en parlais l'autre jour à Mme Cottu qui me disait que son mari en avait fait autant après information prise. Je crois que pour toi dont le nom ne doit être nulle part placé légèrement, il vaut mieux ne pas se presser. Si cela devient une chose bonne et utile nous le saurons bien et il sera toujours tems d'y entrer ; si au contraire cela s'éteint comme je le crois comme une tentative malheureuse, il vaut autant n'avoir pas été au nombre des dupes[99].

Benoît d'Azy n'avait pas tort : la Société française de Statistique universelle, en déclin à partir de 1843, n'a pas survécu à la révolution de 1848[100].

---

# 3. Vers une certaine rigueur

*Laplace et ses héritiers*

Le premier raisonnement de statistique inférentielle de l'histoire date de 1710. Il est dû à John Arbuthnot (1667–1735)[101]. Le résumer en termes modernes nous permettra d'expliciter la terminologie de base de la « preuve statistique ». Arbuthnot part d'une observation : la proportion de garçons parmi les naissances a été supérieure à celle des filles, chaque année pendant quatre-vingt deux ans[102]. Vient alors le raisonnement par l'absurde. Supposons que l'observation soit « due au hasard », c'est-à-dire en l'occurrence que les sexes à la naissance soient aléatoires et indépendants : garçon ou fille avec probabilité 1/2 (c'est l'« hypothèse nulle »). Si cette hypothèse nulle était vraie, alors la probabilité de l'observation serait $2^{-82} = 2,07 \times 10^{-27}$ (c'est la « p-valeur »). Cette probabilité est tellement faible, que l'hypothèse nulle n'est pas acceptable : elle doit être rejetée. Arbuthnot conclut : "From whence it follows, that it is Art, not Chance, that governs".

Arbuthnot raisonnait sur un cas très particulier. La théorie des probabilités, développée par Pierre-Simon Laplace (1749–1827) tout au long de sa carrière, va lui permettre, dès les premiers articles et avant même le théorème central limite, de généraliser le raisonnement en une véritable méthode scientifique[103]. Son premier cas apparaît dans le « Mémoire sur les probabilités », lu le 31 mai 1780[104]. Le sujet est toujours la proportion des sexes à la naissance, qu'il utilisera comme exemple pour ses raisonnements de statistique dans toutes ses publications sur le sujet. Ici, il s'agit de décider au vu d'une observation particulière, si la proportion des garçons est supérieure à 1/2.

> La méthode de l'article précédent donne un moyen fort simple pour obtenir cette probabilité lorsqu'on a un nombre suffisant de naissances ; nous allons l'appliquer à celles qui ont été observées à Paris, et déterminer combien il est probable que les naissances des garçons dans cette grande ville sont plus possibles que celles des filles. Pour cela, nous ferons usage des naissances qui ont eu lieu depuis 1745 jusqu'en 1770, et dont on peut voir la liste dans nos Mémoires pour l'année 1771, page 857. En rassemblant toutes ces naissances, on trouve que, dans l'espace de ces vingt-six années, il est né à Paris 251 527 garçons et 241 945 filles, ce qui donne à très peu près 105/101 pour le rapport des naissances des garçons à celles des filles[105].

Après quelques calculs, Laplace arrive au résultat :

> En repassant des logarithmes aux nombres, on aura, pour la probabilité que $x$ est égal ou moindre que 1/2, une fraction dont le numérateur est peu différent de l'unité et égal à 1,1521, et dont le dénominateur est la septième puissance d'un million ; cette fraction est même un peu trop grande, et, comme elle est d'une petitesse excessive, on peut regarder comme aussi certain qu'aucune autre vérité morale, que la différence observée à Paris entre les naissances des garçons et celles des filles est due à une plus grande possibilité dans la naissance des garçons[106].

En clair, Laplace vient de calculer la p-valeur d'un test de proportion, avec pour hypothèse nulle, que la probabilité de naissance de chaque sexe est 1/2. Il trouve $1,15 \times 10^{-42}$, et la précision impressionnante de son calcul mérite d'être soulignée : le logiciel R donne $1,17 \times 10^{-42}$ ! Bien sûr, on aura remarqué que Laplace à cette époque est encore dans sa « période bayésienne » ; son point de vue évoluera par la suite[107]. Ce qui nous importe ici, c'est la claire perception que Laplace avait d'une démonstration statistique : si la probabilité d'une observation sous l'hypothèse nulle est trop faible, alors l'hypothèse nulle doit être rejetée. Dans un autre mémoire, quelques années plus tard, Laplace montre ce qu'il entend par « petitesse excessive ». Cette fois-ci il s'agit de décider s'il est significatif que la proportion observée à Paris soit inférieure à celle observée à Naples.

---

[101] Stigler, S., 1986, p. 225.
[102] La proportion des sexes à la naissance a fait l'objet de très nombreuse études au cours des siècles : voir Bru, B., 1986, note 53 p. 291, et pour une étude complète, Brian, É. & Jaisson, M., 2007.
[103] Bru, B., 1986, Stigler, S., 1986, Hald, A., 1998, Kuusela, V., 2012.
[104] Laplace, P.-S., 1781.
[105] Laplace, P.-S., 1781, p. 429-430.
[106] Laplace, P. S., 1781, p. 431.
[107] Bru, B., 1986, p. 272.

> On trouve alors la probabilité P, que la possibilité des naissances des garçons à Paris est plus grande que dans le royaume de Naples, égale à 1/100 environ ; il est donc vraisemblable qu'il existe dans ce royaume, comme à Londres, une cause de plus qu'à Paris, qui y facilite les naissances des garçons ; mais la probabilité avec laquelle elle est indiquée par les observations est trop peu considérable encore pour prononcer irrévocablement sur cet objet[108].

Obtenant une p-valeur de 1 %, Laplace considère la conclusion comme vraisemblable, mais n'ose pas pour autant « prononcer irrévocablement ».

Qu'il ait fallu pratiquement un siècle pour que la méthode mise au point par Laplace soit acceptée par les statisticiens, a bien sûr intrigué les historiens. Plusieurs explications ont été avancées. La difficulté technique de la « Théorie analytique des probabilités », qualifiée de « Mont-Blanc mathématique » par de Morgan[109], en est une. La difficulté philosophique à admettre la probabilité laplacienne en est une autre[110]. Dans la période qui nous intéresse, le chef de file des « opposants philosophiques » est Auguste Comte (1798–1857).

> Avant d'examiner la conception fondamentale de Laplace au sujet de l'interprétation cosmogonique des divers caractères généraux que je viens de rappeler, je ne puis m'empêcher de témoigner ici combien tous les bons esprits, étrangers aux préjugés mathématiques, ont dû trouver puérile et déplacée la singulière application du calcul des chances, indiquée d'abord par Daniel Bernouilli, et péniblement complétée ensuite par Laplace lui-même, pour évaluer la probabilité que ces phénomènes ont réellement une cause, comme si notre intelligence avait besoin d'attendre une telle autorisation arithmétique, avant d'entreprendre légitimement d'expliquer un phénomène quelconque bien constaté, lorsqu'elle en aperçoit la possibilité[111].

Dans une longue note de bas de page, Comte insiste : l'application des probabilités à la décision statistique ne peut être que « chimérique ».

> Le calcul des probabilités ne me semble avoir été réellement, pour ses illustres inventeurs, qu'un [pré?]texte commode à d'ingénieux et difficiles problèmes numériques, qui n'en conservent pas moins toute leur valeur abstraite, comme les théories analytiques dont il a été ensuite l'occasion, ou, si l'on veut, l'origine. Quant à la conception philosophique sur laquelle repose une telle doctrine, je la crois radicalement fausse et susceptible de conduire aux plus absurdes conséquences. Je ne parle pas seulement de l'application évidemment illusoire qu'on a souvent tenté d'en faire au prétendu perfectionnement des sciences sociales : ces essais, nécessairement chimériques, seront caractérisés dans la dernière partie de cet ouvrage. C'est la notion fondamentale de la probabilité évaluée, qui me semble directement irrationnelle et même sophistique : je la regarde comme essentiellement impropre à régler notre conduite en aucun cas, si ce n'est tout au plus dans les jeux de hasard. Elle nous amènerait habituellement, dans la pratique, à rejeter, comme numériquement invraisemblables, des événemens qui vont pourtant s'accomplir. On s'y propose le problème insoluble de suppléer à la suspension de jugement, si nécessaire en tant d'occasions. Les applications utiles qui semblent lui être dues, le simple bon sens, dont cette doctrine a souvent faussé les aperçus, les avait toujours clairement indiquées d'avance[112].

Il revient sur le sujet dans le volume 3 de son cours, où il s'en prend cette fois aux applications en statistique médicale, qui avaient été évoquées par Laplace.

---

108 Laplace, P.-S., 1786, p. 325.
109 Bru, B., 1986, p. 283.
110 Bru, B., 1986, p. 280.
111 Comte, A., 1835, p. 286.
112 Comte, A., 1835, p. 287.

> On doit déplorer l'espèce d'encouragement dont les géomètres ont quelquefois honoré une aberration aussi profondément irrationnelle, en faisant de vains et puérils efforts pour déterminer, d'après leur illusoire théorie des chances, le nombre de cas propre à légitimer chacune de ces indications statistiques[113].

Pour autant, la statistique laplacienne ne manquait pas de défenseurs : Bienaymé, Cournot, Fourier, Poisson…[114] : parmi les scientifiques, l'idée qu'une démonstration statistique était non seulement possible mais indispensable, avait été diffusée. Antoine Augustin Cournot (1801–1877) l'exprime très clairement, dès 1828, dans un article où il plaide pour l'introduction de la théorie des probabilités dans l'enseignement.

> Où il n'y a pas de données statistiques les formules du calcul des probabilités sont illusoires ; et partout où ces données mettent en évidence les faits dont nous parlons, les formules trouvent incontestablement leur application[115].

Il exprime les mêmes idées dans son « Exposition de la théorie des chances et des probabilités », plaidant au contraire de Comte pour le « prétendu perfectionnement des sciences sociales ».

> De nos jours, au contraire, la statistique a pris un développement en quelque sorte exubérant ; et l'on n'a plus qu'à se mettre en garde contre les applications prématurées et abusives qui pourraient la décréditer pour un temps, et retarder l'époque si désirable où les données de l'expérience serviront de bases certaines à toutes les théories qui ont pour objet les diverses parties de l'organisation sociale[116].

Si nous avons choisi Cournot de préférence à Bienaymé ou d'autres parmi les héritiers de Laplace, c'est pour un passage de son « Exposition », particulièrement révélateur des précautions qu'il estimait devoir prendre face à Dupin. Ce dernier avait publié l'année précédente, un mémoire sur la proportion des sexes à la naissance, dans lequel comme à son habitude, ses conclusions péremptoires n'étaient pas étayées.

> Les résultats que je viens d'offrir démontrent évidemment qu'il est impossible de considérer comme une quantité *immuable* le rapport entre le nombre des naissances des deux sexes. […] les quinze dernières années, quoique offrant le rapport moyen le plus élevé, présentent une tendance marquée à la diminution successive du rapport des naissances[117].

Or Dupin donnait surtout des rapports bruts, sans les effectifs sur lesquels ils avaient été calculés, ce qui rendait impossible la validation statistique. Le chapitre XIII du livre de Cournot est consacré à l'« Application à des questions concernant les éléments de la population et la durée de la vie ». Il y traite le rapport des sexes à la naissance, en commençant par un long exposé historique, au début duquel il ne manque pas de citer l'article et les données de Dupin, précisant même que « M. Ch. Dupin a eu l'idée heureuse de former deux groupes, l'un des 24 départements maritimes, l'autre des 62 départements de *l'intérieur* »[118]. Ce n'est qu'après une dizaine de pages, que vient l'application proprement dite. Elle concerne la diminution du rapport des sexes, observée pour l'année 1840. Après avoir expliqué le calcul, Cournot conclut :

> Nous aurons $P = 0,1834$, $\Pi = 0,5917$, pour la probabilité que l'écart n'est pas imputable aux anomalies du hasard. Cette probabilité est, comme on le voit, très faible et sans signification possible, quoique l'année 1840 soit, après 1830 et 1828, celle des années de la période où le rapport a atteint sa plus petite valeur[119].

---

113 Comte, A., 1838, p. 419.
114 Bru, B., 1986, p. 280. Voir aussi Bru, B. & M.-F. & Bienaymé, O. 1997, Heyde, C. & Seneta, E. 1977.
115 Cournot, A. A., 1828, p. 7.
116 Cournot, A. A., 1843, p. 181.
117 Dupin, C., 1842, p. 1033.
118 Cournot, A. A., 1843, p. 301.

Ainsi, une des conclusions de Dupin se trouve invalidée, ce sur quoi Cournot se garde bien d'insister !

*Deux amateurs scrupuleux*

À vingt-cinq ans, loin des coteries parisiennes[120], Joseph Honoré Modeste Bigeon (1803–1831) n'a pas la même délicatesse que Cournot : il ose sous-titrer son ouvrage, « où l'on rectifie quelques fausses allégations de M. Ch. Dupin[121] ». Le destin de Bigeon, mort tragiquement à 27 ans « foulé aux pieds des chevaux d'une diligence[122] », n'est pas sans évoquer ceux de ses contemporains Niels Henrik Abel (1802–1829) et Évariste Galois (1811–1832), mais le parallèle tourne court. Malgré une production tout à fait respectable par son volume et son éclectisme[123], Bigeon n'a pas laissé d'empreinte scientifique durable. Lui-même ne s'illusionnait pas sur la portée de son travail statistique, et récusait même le qualificatif de statisticien.

> Vous aurez reçu sans doute avant ma lettre un exemplaire d'un petit ouvrage de peu d'importance à cause de l'incertitude des documents dont j'ai dû faire usage, et que je viens de publier à Toulon. Je vais m'occuper d'un travail du même genre très étendu et, je l'espère fortement, sur toutes les communes du Var. J'ai les tables complètes des décès par âge, par sexes et même par états quant au célibat et au mariage qui m'ont été fournies par l'administration départementale par laquelle j'ai été prié de me charger de cet ouvrage. Il sera peut-être publié par le département, mais c'est fort long et fort ennuyeux et le sentiment de l'utilité dont cela pourrait être m'a seule pu décider à cette entreprise. On me fera, je crois, statisticien, bien que ce soit un genre d'étude bien éloigné de mes goûts. On parle déjà de former une association de 3 ou 4 personnes pour faire la Statistique du département et on me compte parmi, sans quitter mon emploi douanier pourtant[124].

Si le « petit ouvrage de peu d'importance » de Bigeon mérite que l'on s'y arrête, ce n'est pas seulement parce qu'il est un des rares à contester l'autorité de Dupin, c'est surtout à cause de la prudence et de l'honnêteté intellectuelle dont il témoigne. L'objectif est clairement annoncé dès les premières phrases.

> Quelle influence peuvent avoir sur la durée de la vie l'instruction, la richesse, les occupations habituelles des hommes ? Ce sont là des questions d'une assez grande importance et auxquelles il n'est possible de répondre qu'en déterminant pour chaque lieu en particulier la durée moyenne de la vie et en rapprochant sa valeur des autres circonstances locales sous l'influence desquelles les individus considérés sont appelés à vivre.

Pour répondre à la question, Bigeon donne en annexe, d'une part une table de mortalité dont il explique le calcul, d'autre part un tableau comportant huit variables numériques évaluées sur les 86 départements : « vie moyenne », « instruction populaire relative », « revenu foncier moyen », etc. Bigeon est conscient de n'être pas le premier à écrire sur le sujet : il cite bien sûr la théorie analytique des probabilités de Laplace (à trois reprises), mais aussi Rousseau, Godwin, Malthus, Barton, Say, Duvillard, Lalande, Odier, et même Rabelais. Malgré son jeune âge et son isolement,

---

119 Cournot, A. A., 1843, p. 309.
120 Il était employé du service des douanes à Toulon.
121 Bigeon, J. H. M., 1829.
122 Journal des Débats, 3 juin 1831, p. 4.
123 Les éditeurs de Cournot, A. A., 2010, p. 362, jugent peu vraisemblable que le coauteur de Dubois-Aymé pour le « mémoire sur les développées des courbes planes » dont Cournot avait rendu compte, soit aussi l'auteur de l'« Aperçu statistique ». C'est pourtant le cas, comme en atteste la correspondance de Bigeon avec Dubois-Aymé, qui m'a été aimablement communiquée par Pascal Beyls, biographe de Dubois-Aymé (Beyls, P., 2007). Bigeon a également publié deux articles dans les Annales de Physique et de Chimie.
124 Lettre de Bigeon à Dubois-Aymé du 29 juin 1829, transcription communiquée par Pascal Beyls.

Bigeon a donc eu accès à la littérature sur le sujet et l'a comprise. Il n'hésite pas à critiquer ses célèbres devanciers, donnant le plus souvent de bons arguments. Par exemple à propos de Godwin et Malthus :

> Une autre erreur, qui a été commise également par M. Malthus dans son Essai sur le principe de la population, a été d'estimer le nombre des enfans par mariage au moyen de la simple division du nombre des naissances par celui des mariages. Ce calcul est rigoureux dans une population stationnaire, il est suffisant dans une population lentement croissante ; mais lorsque cet accroissement est aussi rapide qu'aux Etats-Unis, il a une énorme influence sur les résultats[125].

C'est un biais numérique du même type que Bigeon reproche à Dupin.

> Il est clair en effet, que si de deux départemens comparés, l'un renferme pour un même nombre d'habitants deux fois plus d'individus qui fréquentent les écoles, mais que dans l'autre la vie moyenne soit deux fois plus longue, il se trouvera à-peu-près dans chacun de ces départemens autant d'individus ayant reçu les bienfaits de l'instruction primaire de sorte que le nombre qui exprimera le rapport d'instruction populaire sera le produit de la vie moyenne multipliée par le nombre d'enfans qui fréquentent les écoles sur un nombre déterminé d'individus ; et si on regarde comme exacts, ce que je suis loin d'admettre, les nombres publiés par M. Dupin, l'instruction relative des départemens sera donnée à-peu-près en divisant par chacun de ces nombres celui qui exprime la durée moyenne de la vie dans le même département. Le changement que cette seule correction apporte dans les nombres de sa carte n'est pas de peu d'importance puisque la vie moyenne varie en France de 29 à 46 ans. En y faisant ces modifications la palme de l'ignorance reste dévolue à la Basse-Bretagne, contrairement aux conclusions du savant académicien[126].

Bigeon est conscient que sa correction ne fournit pas des nombres rigoureusement exacts, et il identifie d'autres sources de biais : la variabilité de la mortalité infantile, de la durée de l'éducation et des structures éducatives. Sa prudence se manifeste dès l'introduction, vis-à-vis des sources.

> Ce travail formé d'après des documens qui n'ont pas toute la certitude désirable, ne peut qu'être incomplet et, sous plus d'un rapport, très-inexact.

Elle manifeste surtout dans ses conclusions (remarquez la déférence ironique de la dernière phrase).

> En considérant, comme l'a fait M. Ch. Dupin dans ses calculs statistiques, la France du Nord d'une part et celle du sud de l'autre, celle-ci comprenant la Bretagne, on trouverait que la vie moyenne est plus longue dans la 1.re partie que dans la seconde, mais je n'oserais en conclure hardiment que c'est l'effet de l'instruction populaire plus répandue au septentrion qu'au midi attendu les exceptions trop formelles que présentent sous le rapport de la durée de la vie les départemens de la Haute-Loire, du Cantal etc. marqués sur sa carte d'une teinte si noire, et j'ajouterai, quoique ce ne soit pas ici le lieu de traiter cette question, au sujet de la liaison admise entre l'instruction et la richesse relatives, que le département de la Seine, malgré sa nuance obscure, est d'un si grand poids dans l'industrie de la France du nord, que je serais, je l'avoue embarrassé de rien décider sur ce point là même, si M. Dupin n'avait pas tranché la question[127].

L'originalité du travail de Bigeon est que ses conclusions se basent sur une méthodologie clairement énoncée. Il est parfaitement conscient de la difficulté d'interpréter des moyennes arithmétiques en cas de valeurs extrêmes, et va même jusqu'à en pondérer certaines.

> Ces résultats sont au reste bien plus incertains que ceux relatifs à la durée de la vie, attendu les variations qu'une seule année produit souvent dans les moyennes de 8 ou 10 années. Je me suis même cru obligé dans le calcul des nombres de la table, lorsqu'une année s'éloignait trop des résultats de toutes les autres, de la faire entrer dans le résultat définitif d'après les règles du calcul

---

125 Bigeon, J. H. M., 1829, p. 49.
126 Bigeon, J. H. M., 1829, p. 11.
127 Bigeon, J. H. M., 1829, p. 9.

des probabilités afin de diminuer proportionnellement son influence[128].

Considérant d'une part les incertitudes sur les données, d'autre part le problème des valeurs extrêmes, Bigeon propose une solution originale pour détecter une liaison entre deux variables. Il détermine les départements dont la valeur est au-dessus de la moyenne pour chacune des deux variables, et considère l'intersection des deux ensembles obtenus : une forte intersection détecte une liaison positive (parler de corrélation serait non seulement anachronique, mais encore partiellement faux). Ses conclusions sont donc basées sur des listes d'intersections, qu'il détaille pages 26 à 28. Par exemple : « Sur les 49 départemens où la vie est longue il en est 23 instruits, 26 ignorans ». Ces affirmations apparemment abruptes, sont en fait rigoureuses car il a pris soin de donner un sens précis à « longue vie », « instruit », etc. Même si la démarche de Bigeon ne va pas jusqu'à un calcul de p-valeur, le simple fait de tenter une réponse méthodologique au problème de l'argumentation statistique, est suffisamment rare à l'époque pour être remarqué. Nous allons voir que d'Angeville a poussé encore plus loin une démarche analogue, qui même selon les critères actuels, peut être considérée comme parfaitement rigoureuse.

Ancien officier de marine, propriétaire terrien, député de l'Ain de 1834 à 1848, auteur de plusieurs rapports sur des sujets divers comme le transport fluvial ou les colonies, Adolphe d'Angeville (1796–1856) n'aurait pas manqué de figurer dans la « statistique morale » de son département, si celle-ci était parue. Il est présenté par les biographes du temps comme :

> un esprit fort disparate et de la famille de ceux auxquels on donne le titre d'*Original*. C'est un mélange du noble et du roturier, de l'homme léger et de l'esprit appliqué; du royaliste et du démocrate[129].

En tout cas quelqu'un qui n'hésite jamais à exercer sa liberté de penser, jusqu'à la tribune de l'assemblée nationale. Son « Essai sur la statistique de la population française[130] » est bien le fruit d'un « esprit appliqué » : 389 pages de texte, dont 8 tableaux comportant un total de 127 variables sur chacun des 86 départements[131], et pas moins de 14 cartes ombrées. Car contrairement à Bigeon, d'Angeville ne s'oppose en rien à Dupin, dont il adopte non seulement le procédé graphique, mais aussi la ligne Saint-Malo Genève : « On serait en effet tenté de croire que deux populations sont venues se heurter en France sur la ligne qui joindrait le *port de St.-Malo* à la ville de *Genève*[132] ». Tout au plus ose-t-il suggérer un échange, qui ferait passer son propre département du bon côté. Il espère probablement que son ouvrage « fait par goût […] et publié à ses propres frais[133] » lui vaudra l'approbation de l'Institut. Il ne s'attendait sans doute pas à la réaction des académiciens au rapport présenté le 5 mars 1838 par Héricart de Thury.

> L'académie des sciences est toujours fort savante, ce qui ne l'empêche pas de s'occuper parfois de graves puérilités. Dans sa séance du 5, elle a entendu l'analyse d'un ouvrage de statistique du comte d'Angerville[sic]. L'auteur y a classé par tableaux toutes les combinaisons que peut présenter la population de la France, en comparant les décès et les naissances par département, et les proportions du développement intellectuel. Enfin, une des principales conclusions de ses immenses recherches a été qu'un pays est d'autant plus éclairé, et que son intelligence est d'autant plus cultivée, qu'il y a plus de fenêtres aux maisons. Messieurs, qui ne s'attendaient pas à ce

---

brusque résultat, n'ont pu réprimer les éclats d'un rire homérique ; ce qui a fort scandalisé le président, M. Becquerel, car il est défendu de rire à l'académie des sciences[134].

Un département est d'autant plus éclairé qu'il y a plus de fenêtres aux maisons ? Oui, c'est bien ce qu'affirme d'Angeville, mais il peut expliquer pourquoi : les ouvertures en tant que base d'imposition sont un indicateur de développement économique, lui-même lié au niveau d'instruction. Les affirmations de d'Angeville, par les raccourcis abrupts auxquels il se livre souvent, sont faciles à caricaturer. Ce qui nous importe ici, c'est qu'elles sont rigoureusement étayées, par une méthodologie que nous allons maintenant examiner.

Comme Bigeon, d'Angeville est très conscient d'une part de l'imprécision de ses chiffres, d'autre part de la sensibilité des moyennes aux valeurs extrêmes. S'il ajoute un « département moyen » à toutes ses variables, c'est sans doute pour sacrifier à la mode du temps (l'homme moyen de Quetelet date de 1831[135]) ; mais il exprime ses doutes[136], et ne l'utilise jamais dans son argumentation. Par contre, toutes les variables utiles sont systématiquement couplées dans ses tableaux à leurs statistiques de rang.

> Pour les calculs qui ont quelque importance, nous avons donné à tous les départements un numéro d'ordre relatif […] ; de cette manière, on peut voir du premier coup d'œil l'ordre des divers départements entr'eux, pour les faits qui s'appliquent aux recherches les plus utiles.

Écartant la Corse « à raison de sa situation toute exceptionnelle », il reste 85 départements,

> […] dont le cinquième est 17. C'est ce dernier nombre qui compose une série. Ainsi lorsque nous parlons de la série où la population est la plus dense, nous voulons parler des 17 départements où il y a le plus de population par myriamètre carré. Cette explication était nécessaire pour l'intelligence de la suite du travail ; elle s'applique à tous nos calculs[137].

Un des avantages de cette notion de série, est de donner une base plus rigoureuse aux cartes ombrées : les siennes utilisent cinq teintes, et les départements y reçoivent la couleur de la série à laquelle ils appartiennent[138]. Plus important, les premières et cinquièmes séries fondent l'argumentation statistique de d'Angeville. Prenons l'exemple des fenêtres, qui a tant diverti l'académie. D'Angeville dit : « Sur les 17 départements de la série où il y a le plus d'ignorance, 12 sont aussi de celle où il y a le moins d'ouvertures imposables aux maisons ». Implicitement, il considère donc qu'une intersection de 12 départements entre deux séries de 17, est significativement grande. Mais si les deux séries sont indépendantes l'une de l'autre (hypothèse nulle), la taille de leur l'intersection doit suivre une loi géométrique de paramètres 17, 68, 17. On peut donc calculer exactement la probabilité qu'elle soit supérieure ou égale à 12, c'est-à-dire la p-valeur. On trouve $2.08 \times 10^{-7}$, donc d'Angeville a raison : une intersection de 12 départements entre deux séries de 17 est significative. La procédure ci-dessus, d'usage très courant de nos jours, porte le nom de « test exact de Fisher[139] ». Ce n'est certainement pas le test d'indépendance le plus puissant, mais il est parfaitement correct. Il n'a été formulé ainsi qu'un siècle après d'Angeville ; néanmoins le calcul, qui ne fait intervenir que des coefficients combinatoires, est plus simple que l'application de la théorie de Laplace, et tout à fait à la portée ne n'importe quel probabiliste de l'époque. Il est possible de contrôler chaque affirmation de d'Angeville, en calculant sa p-valeur pour le test exact de Fisher ; à une exception près (citation suivante), toutes ces p-valeurs sont inférieures à 5 %, preuve que l'intuition de d'Angeville sur ce qu'est une intersection significative est fiable. Voici un exemple, dans lequel d'Angeville ne cache pas son propre étonnement à voir un préjugé très répandu infirmé

---

134 Ollivier, J., 1838, p. 228.
135 Stigler, S. 1999, p. 59.
136 D'Angeville, A., 1837, p. 14.
137 D'Angeville, A., 1837, p. 20.
138 D'Angeville, A., 1837, p.349.
139 Agresti, A., 1992, p. 134.

par ses données. Les p-valeurs sont rajoutées entre crochets.

> Voici ce que les faits indiquent : sur les 17 départements où il y a le plus d'ignorance, 7 sont de la série où il y a le moins d'accusés de crimes [P = 0,022], tandis qu'un seul, celui des *Pyrénées-Orientales* figure dans celle où l'on compte le plus de ces accusés [P = 0,092]. […] Dès lors, il nous a bien fallu reconnaître une vérité : c'est que la criminalité n'est en aucune manière déterminée par le défaut d'instruction. […] Quel n'a pas été notre étonnement, lorsque nous avons vu que les 32 départements de la France du Nord, qui sont si éclairés, contiennent 13 des 17 départements de la série qui présente le plus d'accusés de crimes, tandis que le midi, c'est-à-dire 53 départements, n'en renferme que 4 [P = 3,6x10$^{-4}$][140].

Un autre exemple illustre bien l'attitude scientifique de d'Angeville.

> Si l'on examine tous les départements baignés par la mer, on voit que pas un n'est de la série des 17 départements où il y a le plus de goitres [P = 0.0043]. […] Y aurait-il dans le voisinage de la mer quelque chose qui s'oppose à ce mal[141] ?

La découverte par Coindet de l'effet bénéfique de l'iode sur les goitres commençait à être connue[142], mais d'Angeville l'ignorait. Pourtant, il avait bien détecté l'effet de la mer, par sa seule méthode statistique. Sans faire de d'Angeville le nouveau père de la statistique inférentielle, on peut au moins le compter parmi les précurseurs de la statistique non paramétrique[143], pour son utilisation systématique des statistiques de rang, et l'usage rigoureux qu'il en fait.

## 4. L'extension du phénomène

En fin observateur de la société de son temps, Balzac ne pouvait pas s'être trompé : il y a bien eu une mode de la statistique en France. De l'examen de la presse et de la littérature de l'époque, il semble possible de conclure qu'elle a démarré en 1827, à la suite de la carte ombrée de Dupin. Il est beaucoup plus difficile d'en préciser l'extension. Nous poserons la question en termes d'étendue géographique, puis sociologique, puis historique.

Comme l'attestent le témoignage de Dupin et l'épisode de la statistique morale d'Andraud, le phénomène a pu démarrer au sud de la ligne Saint-Malo Genève, en réaction aux accusations d'obscurantisme. Cinq sociétés savantes comportant le mot statistique dans leur intitulé (elles ne s'occupaient pas que de statistique) ont été créées en province entre 1827 et 1838 : à Marseille en 1827, Bourges en 1834, Niort en 1836, Valence en 1837, Grenoble en 1838, donc toutes au sud de la ligne. En revanche, le succès parisien des « stationnaires », pièce écrite en septembre 1827, montre que le retentissement était national. Il y a peut-être eu plus d'articles de statistique locale dans le midi, mais dès 1828, de nombreuses parutions nationales comportent le mot statistique dans leur titre, probablement pour des raisons publicitaires. Il semble donc qu'au moins la vogue du mot ait bien été nationale. Il serait intéressant d'en examiner le retentissement international, plus particulièrement en Belgique et en Angleterre. On sait que l'organisation de la discipline doit beaucoup aux réunions qui ont eu lieu à Londres en 1833, à l'occasion de la visite de Quetelet[144]. Mais Babbage et Quetelet avaient fait connaissance à Paris en 1826[145], et avaient régulièrement correspondu depuis lors. Nous ignorons si le phénomène dont il est question ici, largement indépendant de la statistique institutionnelle, a pu avoir un écho chez nos voisins.

Combien de Magis et d'Homais y a-t-il eu en France ? Le personnage, s'il est bien identifié

---

140 D'Angeville, A., 1837, p. 69.
141 D'Angeville, A., 1837, p. 55.
142 Mérat, F.-V. & de Lens, A. J., 1831, p. 632.
143 Droesbeke J.-J. & Fine, J. 1996, chapitre I.
144 Hacking, I, 1990, p. 61, Babbage, 1864, p. 434, Quetelet, A., 1873, p. 156.
145 Quetelet, A., 1873, p. 152.

dans la littérature, y reste relativement rare. Presque tous les romanciers ont utilisé le mot statistique, mais seuls Balzac, Flaubert et Zola ont peint des statisticiens. Il n'y a chez Balzac qu'un personnage annonçant son intention de publier un ouvrage de statistique, sur les quelque 2000 qu'il a créés pour la Comédie Humaine[146]. On trouve bien Homais chez Flaubert, mais Bouvard et Pécuchet, ces boulimiques de connaissance, n'évoquent la statistique que marginalement[147]. Outre les « six membres de la société de statistique », on n'en compte que deux qui aient un nom, sur les quelque 1200 personnages des Rougon-Macquart[148]. La question peut aussi être abordée au moyen du recensement des sociétés savantes. N'en déplaise à Labiche, il n'y a jamais eu de société de statistique à Vierzon. On compte lors du recensement de 1846, 134 sociétés en province, dont seulement 4 sont identifiées comme « sociétés de statistique[149] » (celle de Bourges avait disparu dès 1840) ; pour comparaison, il y avait 22 sociétés médicales ou pharmaceutiques. On doit cependant observer que de nombreuses autres sociétés publiaient des statistiques sans afficher l'activité dans leur titre, ni même dans leurs objectifs. Néanmoins, il ne semble pas que la production d'ouvrages intitulés « statistique » ait dépassé quelques pour cent du total. Par exemple, à la Société des sciences, Belles-Lettres et Arts du département du Var, dont il est membre depuis 1828, Bigeon est le seul à pouvoir être identifié comme statisticien[150]. Les quatre sociétés de statistique de province totalisent 242 membres en 1846. Nous avons vu ce qui pouvait se dire des chiffres avancés par César Moreau à Paris. Il ne semble pas que le nombre d'érudits de province ayant publié au moins un ouvrage de statistique ait pu dépasser un ou deux milliers.

Jusqu'à quand la mode de la statistique a-t-elle duré ? Dès 1827, la statistique est incarnée pour le grand public par Dupin, au point que les journaux considèrent souvent qu'il en est l'initiateur. Ainsi en 1829 le Pirate parle de « cette science toute nouvelle […] que M. Charles Dupin a inventée à son profit[151] ». Près de quarante ans plus tard, le Tintamarre considère que « M. le baron Charles Dupin est le créateur de cette science folâtre que quelques esprits sérieux se sont plûs à qualifier de statistique à leurs moments perdus[152] ». Passés le décès du « créateur » en 1873 et les hommages consécutifs, le nom statistique n'est plus associé au sien dans la presse. Pourtant en 1880, « vivre et mourir en statistiquant[153] », faisait encore rire dans une opérette à succès. Le Vabre de Zola date de 1882[154], et il semble que ce soit le dernier statisticien identifié en tant que personnage comique. Cela ne signifie pas que l'on ait cessé de parler de statistique, bien au contraire : plus de 56 % des numéros du Journal des Débats utilisent le mot entre 1880 et 1900, contre 25 % entre 1828 et 1848. Vers 1896, l'affaire Dreyfus révèle au public une acception qui ne prête plus à sourire : c'est de la section de statistique du ministère des armées qu'est venue l'accusation, et comme on le saura plus tard, la fabrication des pièces[155]. Il semble qu'alors, la période d'exubérance soit passée : le mot statistique est désormais employé comme il l'est de nos jours, c'est-à-dire le plus souvent bien loin de ce que les scientifiques en comprennent. Que cela ne nous empêche pas, aujourd'hui comme hier… d'en rire !

> Mais, comme dans tout pays civilisé il faut, pour la bonne tenue des statistiques, qu'il y ait une moyenne de cocus déterminée, assez imposante pour que nous ne soyons pas dans une situation inférieure vis-à-vis des autres nations, je vous conseillerais de compléter votre proposition en disant que, pour parfaire les manquants, le titre de cocu serait attribué d'office à tous les membres de l'Académie des sciences morales et politiques, et, si cela ne suffit pas, à tous les membres de la Société d'économie politique[156].

---

146 Cerfberr, A. & Christophe, J., 1887.
147 Flaubert, G., 1881, p. 171.
148 Ramond, F. C., 1901, préface.
149 Comte, A.-J., 1846, p. 1019.
150 Compte rendu des travaux de la Société des sciences Belles-Lettres et Arts du département du Var séant à Toulon, pendant les années 1830 et 1831, Toulon, imprimerie de Baume 1832.
151 Le Pirate, 6 septembre 1829.
152 Le Tintamarre, 5 mai 1868.
153 Offenbach, J., et al., 1880, p. 38.
154 Zola, É, 1882, p. 96.
155 Journal des Débats, 15 septembre 1896, p. 2. Voir Laurent, S., 2010 pour l'histoire des relations entre statistique et renseignement en France.

## Bibliographie

---

156 Allais, A., 1895, p. 185.